\documentclass{amsart}
\usepackage[all]{xy}
\usepackage{graphicx}
\usepackage{psfrag}
\usepackage{amsthm}
\usepackage{amscd}
\usepackage{amsfonts}
\usepackage{amstext}
\usepackage{amssymb}
\usepackage{amsmath}
\usepackage{amsxtra}
\usepackage{plain}
\usepackage[all]{xy}

\newtheorem{prop}{\bf Proposition}[section]

\newtheorem{thm}[prop]{{\bf Theorem}}

\numberwithin{equation}{section}

\newenvironment{pf}{{\bf proof }}{\qed\endtrivlist}

\newcommand{\Q}{\mathbb{Q} }
\newcommand{\Z}{\mathbb{Z} }

\newcommand{\C}{\mathbb{C} }
\newcommand{\R}{\mathbb{R} }

\newcommand{\ind}{\operatorname{ind}}

\newcommand{\Ch}{\operatorname{Ch}}

\newcommand{\Hom}{\operatorname{Hom}}

\begin{document}

\title{An appendix to a paper by B. Hanke and T. Schick}
\author{Mostafa ESFAHANI ZADEH}
\address{Mostafa Esfahani Zadeh.\newline
 Mathematisches Institute,
Georg-August-Universit\"{a}t,\newline
G\"{o}ttingen, Germany\and \newline 
Institute for Advanced Studies in Basic Sciences,\newline
Zanjan, Iran}
\email{zadeh@uni-math.gwdg.de}
\maketitle
\begin{abstract}
In this short note we apply methods introduced by B. Hanke and T. Shick to prove the 
vanishing of (low dimensional-)higher $A$-genera for spin manifolds admitting a positive 
scalar curvature metric. Our aim is to provide a short and unified proof for this beautiful 
result without using the strong Novikov conjecture. 
\end{abstract}

\section{Introduction}
Let $G$ be a discrete group and $BG$ its classifying space. The strong Novikov conjecture 
asserts the injectivity of the rationalized assembly map 
\[A\colon K_*(BG)\otimes\Q\to K_*(C^*_{max}G)\otimes Q~.\] 
 Recently B. Hanke and
T. Schick have proved the strong Novikov conjecture for those class in $K_*(BG)$ belonging 
to the dual of $\Lambda^*(G)\subset H^2(BG;\Q)$. Here $\Lambda^*(G)$ denote the subring of
$H\sp*(BG;\Q)$ generated by elements in $H^2(BG;\Q)$ and the duality is defined by 
means of the homological Chern character $Ch:K_*(BG)\to H_*(BG;\Q)$, 
c.f. \cite[section 11]{BaDo}. 
Their proof is based on several sophisticated techniques in higher 
index theory and in particular a fundamental role is 
played be an assembling techniques that they have invented in \cite{HaSc1, HaSc2}. 
In this short note we are interested in an application of the strong Novikov conjecture. 
Let $M$ be a closed spin manifold 
and let $f\colon M\to BG$ be any continuous map. It is proved in \cite[theorem 3.5]{J.Ros1}
that if the strong Novikov conjecture holds for $G$ and if $M$ 
admit a metric with everywhere positive scalar curvature
then for all $a\in H^*(BG,\Q)$,
\begin{equation}\label{novvan}
 \langle\hat A(TM)\cup f^*(a),M\rangle=0~.
\end{equation}
As an immediate corollary of this result, if $f_*([M])\neq0$ in $H_*(BG,\Q)$
then $M$ does not admit a metric with everywhere positive scalar curvature.
For $G$ hyperbolic this corollary was proved by M. Gromov and B. Lawson by means 
of their relative index theorem (see \cite[theorem 13.8]{GrLa3}). In this short note we 
follow \cite{HaSc3} to give a 
proof for vanishing relation \eqref{novvan} provided $a\in \Lambda^*(G)$. We do not deal with the 
strong Novikov conjecture and do not use the assembling construction.  
This makes all arguments easier 
to follow and provides a short proof for the vanishing higher genera in low dimension. 
As the title suggests this note should be considered as an appendix to the paper \cite{HaSc3}. 
The author has merely applied the method of this paper to prove a result which is over looked in 
\cite{HaSc3}. 

\section{The vanishing theorem}
\begin{thm}
Let $\Lambda^*(G)$ denote the subring of $H\sp*(BG;\Q)$ generated
by elements in $H^2(BG;\Q)$.
Let $M$ be a closed spin manifold admitting a riemannian metric with
positive scalar curvature. Then for all  $a$ in $\Lambda^*(G)$ 
\[\langle\hat A(TM)\cup f\sp*(a),[M]\rangle=0~.\]
In particular if $f_*[M]$ belongs to the dual of $\Lambda^*(G)$
then $M$ cannot carry a metric with positive scalar curvature.
\end{thm}
\begin{pf}
By multiplying with an appropriate integer, we may assume $a\in H^2(M,\Z)$.
Let $\ell\to BG$ be the smooth unitary complex line bundle classifying by $a$. The unitary bundle 
$L:=f^*\ell$ is then classified by 
$f^*(a)$, i.e. the first Chern class $C_1(L)$ equals $f^*(a)$. Fix a unitary
connection on $L$ and denote its curvature by $\omega\in\Omega^2(M)$;
a closed differential $2$-form which represents $C_1(L)$. 
Denote by $\tilde M$ the covering 
$\tilde M\stackrel{\pi}{\rightarrow}M$ which is classified by $f$. 
The discrete group $\Gamma:=\pi_1(M)/\ker f_*$ is the Deck transformation 
group of $\tilde M$ and $\tilde M/\Gamma=M$. 
Let $\tilde L:=\pi^*(L)$ denote the lifting of $L$ to $\tilde M$. 
Since the universal cover of $BG$
is contractible, and by naturality, the bundle $\tilde L$ is trivial. 
Fix a unitary isomorphism $\tilde L\simeq\tilde M\times \C$. With respect to
this isomorphism, the lifted connection on $\tilde L$ takes the form $d+\eta$
where $\eta\in\Omega^1(\tilde M)$. Since the structural group of $L$
is abelian, the curvature of this connection is $d\eta$. 
Though $\eta$ is not in general $\Gamma$-invariant its curvature is and
by naturality $d\eta=\pi^*(\omega)$. \newline
For $0\leq t\leq1$ consider the family $d+t\eta$ of unitary
connections on $\tilde L$. The curvature of this family is given
by $\tilde\omega_t:=t\pi^*(\omega)$ and satisfies the following inequality
\begin{equation}\label{avnam}
\|\tilde\omega_t\|\leq c.t~.
\end{equation}
Here $c$ is a constant depending on the geometry of the vector bundle
$L$. 
Consider the vector bundle 
\[\mu:=\tilde M\times_\Gamma\ell^2(\Gamma)~.\]
The connection $d+t\eta$ induces a connection $\nabla_t$ on this bundle.with curvature 
$\omega_t$. 
In fact a smooth section of $\mu$ may be considered as a $L^2$-smooth
section of $\tilde M\times \C$ and 
the action of $\nabla_t$ on such a section coincides with the action of $d+t\eta$. 
It is clear from 
this description that the connection $\nabla_t$ is unitary with respect to the inner
product of $\ell^2(\Gamma)$ and its curvature $\omega_t$ satisfies the following inequality
\begin{equation}\label{dovnam}
\|\omega_t\|\leq c.t~,
\end{equation}
where $c$ is the constant of \eqref{avnam}. In \cite{HaSc3} the bundle $\mu$ and the connection
$\nabla_t$ are used to construct a Hilbert $A_t$-module bundle $V_t$. We describe briefly 
this construction in below. Fix a point $x\in M$ and a $\tilde x\in\tilde M$. So one can identify 
naturally $\mu_x$ with $\ell^2(\Gamma)$. Let $y\in M$
be an arbitrary point and $\gamma$  be a piecewise smooth curve from $x$ to $y$. Let
$A_t(\gamma)\in\Hom(\mu(x),\mu(y))$ denote the parallel translation
with respect to $\nabla_t$. For a fix $y$, denote by $V_t(y)$ the linear-norm
completion of the set of all such $A_t(\gamma)$'s. These linear spaces
form a smooth vector bundle $V_t$. Notice that $A_t:=V_t(x)$
is in fact a $C^*$-algebra
and has a right module action on each fiber of $V_t$ given by precomposition.
Moreover the following pairing
\begin{equation}\label{inerpro}
\langle A_t(\gamma),A_t(\gamma')\rangle:=A_t(\gamma)^{-1}\circ A_t(\gamma')
\end{equation}
defines an $A_t$-valued bilinear form which supplies $V_t$ with the structure of a Hilbert
$A_t$-module bundle. Let $v$ be an element in $\mu_{y}$ which is the parallel 
translation of $u\in V_{x}$ along 
the curve $\gamma$ from $x$ to $y$ and let $\beta$ be a piecewise smooth curve from 
$y$ to $z$. Parallel translation of $u$ along $\beta\circ\gamma$ define $v'\in V_{t,z}$. 
The family $P(\beta)$ given by 
\begin{equation}
P(\beta):V_{t,y}\to V_{t,z}~;\hspace{5mm}P(\beta)(A(\gamma))(v)=v'
\end{equation}
satisfy the properties of parallel translation. Clearly  
these parallel translations are $A_t$-linear and are unitary with respect to \eqref{inerpro}. 
Therefore they define 
an $A$-linear unitary connection on $V_t$ whose curvature
$\Omega_t$ is given by  
\begin{equation}\label{socur}
\Omega_t(y)(A)=A\circ \omega_t(y)~,\hspace{5mm}A\in V_{t,y}
\end{equation}
and satisfies again the following inequality
\begin{equation}\label{somnam}
\|\Omega_t\|\leq c.t~,
\end{equation}
The algebra $A_t$ has a complex-valued natural
trace $\tau_t$ given by the
following formula, c.f. \cite[Lemma 2.2]{HaSc3}
\[\tau_t(A_t(\gamma))=\langle A_t(\gamma)(1_e),1_e\rangle~.\]
Notice that the identification $\mu_x=\ell^2(\Gamma)$ is used here. 
Let $D^{V_t}$ denote the spin Dirac operator on $M$ twisted by the bundle $V_t$. The
Mishchenko-Fomenko index $\ind D^{V_t}$ takes its value in $K_0(A_t)$ which is
acted on by $\tau_t$. The following index formula gives a topological expression 
for $\tau_t(\ind D^{V_t})$ (see \cite[Theorem 6.9]{Sc-L^2})
\[\tau_t(\ind D^{V_t})=\langle\hat A(TM)\cup\Ch_{\tau_t}(V_t),[M]\rangle~,\]
where $\Ch_{\tau_t}(V_t)=\tau_t(\exp(\Omega_t))$. 
From \eqref{socur} one has $\tau_t(\Omega_t)=\tau_t(t\pi^*f^*(a))=tf^*(a)$, so  
\[\Ch_{\tau_t}(V_t)=1+tf^*(a)+t^2f^*(a)\wedge f^*(a)+\dots\]
is a polynomial function of $t$ of degree at most $m=\dim M$ with coefficient in $\Lambda^*T^*M$. 
Consequently 
\[\tau_t(\ind D^{V_t})=\langle\hat A(TM),[M]\rangle+
t\langle\hat A(TM)\cup f^*(a),[M]\rangle+\dots\in\R[t].\]
Inequality \eqref{somnam} and the Lichnerowicz formula imply the
vanishing of $\ind D^{V_t}$ for small values of $t$. The vanishing of above polynomial for small values 
of $t$ implies the vanishing of all its coefficients. In particular  
\[\langle\hat A(TM)\cup f^*(a),[M]\rangle=0,\]
which is the desired relation and completes the proof of the theorem for $a\in H^2(M,\Z)$. Now let 
$a=a_1\cup\dots\cup a_k$ where each $a_j$ is an element in $H^2(BG,\Q)$. Consider
the line bundle $\mathcal L:=L_1\otimes\dots\otimes L_k$ where $L_j$ is the line bundle
classified by $f^*(a_j)$. Proceeding as in above by using the line bundle $\mathcal L$ 
(instead of $L$) and 
the connection $d+t_1\eta_1+\dots+t_k\eta_k$ (instead of $d+t\eta$), one deduces the 
vanishing of a polynomial function of variables $t_j$ for $1\leq j\leq k$.
The  coefficient of $t_1t_2\dots t_k$ in this polynomial is
\[\langle\hat A(TM)f^*(a_1\cup\dots\cup a_k),[M]\rangle.\]
Therefore the above expression must be zero. This complete the proof of the theorem.
\end{pf}

\end{document}